\documentclass[12pt,leqno]{amsart}

\usepackage[dvips]{graphicx}

\usepackage{fullpage}
\usepackage{amsmath,amssymb,amsthm}
\usepackage{enumerate}
\usepackage{float}

\newtheorem{thm}{Theorem}
\newtheorem{prop}[thm]{Proposition}
\newtheorem*{RP}{Piecewise $D$-finite Reconstruction Problem}

\DeclareMathOperator{\dd}{d}
 
\DeclareMathOperator{\Op}{\mathfrak{D}} 

\DeclareMathOperator{\diag}{diag}

\newcommand{\der}[2]{#1^{(#2)}}
\newcommand{\nullsp}[1]{\ensuremath{\mathcal{N}_{#1}}} 
\newcommand{\f}[1]{\mbox{$#1$}}
\newcommand{\fun}{\ensuremath{f}}
\newcommand{\isdef}{\ensuremath{\stackrel{\text{def}}{=}}}
\newcommand{\reals}{\ensuremath{\mathbb{R}}}

\newcommand{\fn}{\ensuremath{g}} 
\newcommand{\theorder}{\ensuremath{r}} 
\newcommand{\polymaxorder}{\ensuremath{k}} 
\newcommand{\jc}{\ensuremath{a}}
\newcommand{\meas}[1][k]{\ensuremath{m_{#1}}}
\newcommand{\fwm}{\ensuremath{\mathcal{P}}} 

\newcommand{\jp}{\ensuremath{x}}
\newcommand{\np}{\ensuremath{N}} 
\newcommand{\polycoeff}{\ensuremath{b}} 

\newcommand{\nn}[1]{\widetilde{#1}} 

\newcommand{\nc}{\newcommand}
\nc{\be}{\begin{equation}}
\nc{\ee}{\end{equation}}
\nc{\Rr}{{\mathbb R}}
\nc{\Nn}{{\mathbb N}}
\nc{\Ff}{{\mathcal F}}
\def\phi{\varphi}

\begin{document}

\title{
An ``algebraic'' reconstruction of piecewise-smooth functions from integral measurements}
\thanks{Supported by ISF and the Minerva Foundation}

\author{D.Batenkov}
\address{Department of Mathematics, Weizmann Institute of Science, Rehovot 76100, Israel}
\email{dima.batenkov@weizmann.ac.il}

\author{N.Sarig}
\email{niv.sarig@weizmann.ac.il}

\author{Y.Yomdin}
\email{yosef.yomdin@weizmann.ac.il}


\maketitle

%
%

\begin{abstract}
This paper presents some results on a well-known problem in
Algebraic Signal Sampling and in other areas of applied
mathematics: reconstruction of piecewise-smooth functions from
their integral measurements (like moments, Fourier coefficients,
Radon transform, etc.). Our results concern reconstruction (from
the moments or Fourier coefficients) of signals in two specific classes:
linear combinations of shifts of a given function, and ``piecewise
$D$-finite functions" which satisfy on each continuity interval a
linear differential equation with polynomial coefficients. In each
case the problem is reduced to a solution of a certain type of
non-linear algebraic system of equations (``Prony-type system"). We
recall some known methods for explicitly solving such systems in
one variable, and provide extensions to some multi-dimensional cases.
Finally, we investigate the local stability of solving the
Prony-type systems.
\end{abstract}

\section{Introduction}

It is well known that the error in the best approximation of a
$C^k$-function $f$ by an $N$-th degree Fourier polynomial is of
order ${C\over {N^k}}$.
The same holds for algebraic polynomial
approximation and for other basic approximation tools (see e.g. \cite[Chapters IV, VI]{Natanson1949}, and \cite[Vol.I, Chapter 3, Theorem 13.6]{zygmund1959trigonometric}). However,
for $f$ with singularities, in particular, with discontinuities,
the error is much larger: its order is ${C\over {\sqrt N}}$.
Considering the so-called Kolmogorow $N$-width of families of
signals with moving discontinuities one can show that {\it any
linear approximation method provides the same order of error, if
we do not fix a priori the discontinuities' position} (see
\cite{Ett.Sar.Yom}, Theorem 2.10). Another manifestation of the
same problem is the {\it ``Gibbs effect" - a relatively strong
oscillation of the approximating function near the
discontinuities}. Practically important signals usually do have
discontinuities, so the above feature of linear representation
methods presents a serious problem in signal reconstruction. In
particular, it visibly appears near the edges of images compressed
by JPEG, as well as in the noise and low resolution of the CT and
MRI images.

Recent non-linear reconstruction methods, in particular,
Compressed Sensing (\cite{Can,Don}) and Algebraic Sampling
(\cite{Vet1,Vet2,Vet3,Mil,Put}), address this problem in many
cases. Both approaches appeal to an a priori information on the
character of the signals to be reconstructed, assuming their
``simplicity" in one or another sense. Compressed sensing assumes
only a sparse representation in a certain (wavelets) basis, and
thus it presents a rather general and ``universal" approach.
Algebraic Sampling usually requires more specific a priori
assumptions on the structure of the signals, but it promises a
better reconstruction accuracy. In fact, we believe that
ultimately the Algebraic Sampling approach has a potential to
reconstruct ``simple signals with singularities" as good as smooth
ones. The most difficult problem in this approach seems to be estimating
the accuracy of solution of the nonlinear systems arising.

Our purpose in this paper is to further substantiate the Algebraic
Sampling approach. On one hand, we present two algebraic reconstruction
methods for generic classes of signals. The first one is reconstruction
of combinations of shifts of a given function from the moments and
Fourier coefficients (Section \ref{sec:shifts}). The second one concerns
piecewise $D$-finite moment inversion (Section \ref{sec:dfinite}).
On the other hand, we consider some typical nonlinear systems arising in these reconstruction schemes. We describe the methods of their solution (Section \ref{sec:Prony.System}), and provide some explicit
bounds on their local stability (Section \ref{sec:prony}). We also present results of some numerical experiments in Section \ref{sec:numerical}.

Our ultimate goal may be stated in terms of the following conjecture
(which seems to be supported also by the results of \cite{Eck,Kve,Tad1,Tad2,Vet3}):

{\it There is a non-linear algebraic procedure reconstructing any
signal in a class of piecewise $C^k$-functions (of one or several
variables) from its first $N$ Fourier coefficients, with the
overall accuracy of order ${C\over {N^k}}$. This includes the
discontinuities' positions, as well as the smooth pieces over the
continuity domains.}

Recently \cite{BatYomFourer} we have shown that at least half the conjectured
accuracy can be achieved. However, the question of maximal possible accuracy
remains open. Our results presented in this paper can be considered as an
additional step in this direction.

\section{Linear combinations of shifts of a given function}\label{sec:shifts}

Reconstruction of this class of signals from sampling has been
described in \cite{Vet1,Vet2}. We study a rather similar problem
of reconstruction from the moments. Our method is based on the
following approach: we construct convolution kernels dual to the
monomials. Applying these kernels, we get a Prony-type system of
equations on the shifts and amplitudes.

Let us restate a general reconstruction problem, as it appears in
our specific setting. We want to reconstruct signals of the form
\be\label{eq:the_model} F(x)=\sum_{j=1}^Na_j f(x-x^j) \ee where $f$ is a
known function of $x=(x_1,\ldots,x_d)\in {\mathbb R}^d$, while the
number $N$ of the shifts and the form (\ref{eq:the_model}) of the
signal are known a priori. The parameters $a_j, \
x^j=(x^j_1,\dots,x^j_d), \ j=1,\dots N$ are to be found from a finite
number of ``measurements", i.e. of linear (usually integral) functionals
like polynomial moments, Fourier moments, shifted kernels,
evaluation over some grid etc.

In this paper we consider only linear combinations of shifts of
one known function $f$. Reconstruction of shifts of several functions
based on ``decoupling" via sampling at zeroes of the Fourier transforms
of the shifted functions is presented in \cite{Sar.Yom3}.

\medskip

In what follows $x=(x_1,\dots,x_d), t=(t_1,\ldots,t_d) \in {\mathbb R}^d$,
\ $j$ \ is a scalar index, while $k=(k_1,\dots,k_d), \ i=(i_1,\dots,i_d)$
and $n=(n_1,\dots,n_d)$ are multi-indices. Partial ordering of
multi-indices is given by $k \leq k' \Leftrightarrow k_p \leq
k'_p, \ p=1,\dots,d.$

\medskip

Assume that a multi-sequence of functions
$\phi= \{\phi_k(t)\},\; t\in {\mathbb R}^d, \; k \geq (0,\dots,0)$
is fixed. We consider the measurements $\mu_k(F)$ provided by

\be\label{eq:measurments}\mu_k(F)=\int F(t)
\phi_k(t)\dd t, \ k \geq (0,\dots,0).\ee Our approach now works as follows:
given $f$ and $\phi=\{\phi_k(t)\}$ we now try to find an ``$f$-convolution
dual" system of functions $\psi = \psi_k(t)$ in a form of certain
``triangular" linear combinations

\be\label{eq:triangular_system} \psi_k(t)=\sum_{0\leq i\leq k}
C_{i,k}\phi_i(t), \ k\geq (0,\dots,0). \ee More accurately, we try
to find the coefficients $C_{i,k}$ in (\ref{eq:triangular_system})
in such a way that \be\label{eq:dual_convolution} \int f(t-x)\psi_k(t)\dd t =
\phi_k(x).\ee We shall call a sequence $\psi = \{\psi_k(t)\}$
satisfying (\ref{eq:triangular_system}), \
(\ref{eq:dual_convolution}) \ $f$ - convolution dual to $\phi$.
Below we shall find explicitly convolution dual systems to the usual
and exponential monomials.

\medskip

We consider a general problem of finding convolution dual
sequences to a given sequence of measurements as an important step
in the reconstruction problem. Notice that it can be generalized
by dropping the requirement of a specific representation
(\ref{eq:triangular_system}): $\psi_k(t)=\sum_{i=0}^k
C_{i,k}\phi_i(t)$. Instead we can require only that $\int
f(t)\psi_k(t) \dd t$ be expressible in terms of the measurements
sequence $\mu_k$. Also $\phi_k$ in (\ref{eq:dual_convolution}) can
be replaced by another a priori chosen sequence $\eta_k$. This
problem leads, in particular, to certain functional equations,
satisfied by polynomials and exponential functions, as well as exponential
polynomials and some kinds of elliptic functions (see \cite{Sar.Yom2}).

Now we have the following result:

\begin{thm}\label{thm:conv_dual} Let a sequence $\psi = \psi_k(t)$ be
$f$-convolution dual to $\phi$. Define $M_k$ by $M_k=\sum_{0\leq i\leq k}
C_{i,k}\mu_i.$ Then the parameters $a_j$ and $x^j$ in
(\ref{eq:the_model}) satisfy the following system of
equations (``generalized Prony system"):
\be\label{eq:gnrlzd_prony} \sum_{j=1}^N a_j\phi_k(x^j)=M_k, \ \
k\geq (0,\dots,0).\ee
\begin{proof}
By definition of $M_k$ and via (\ref{eq:triangular_system})
and (\ref{eq:dual_convolution}) we have for each $k\geq (0,\dots,0)$

$$M_k=\sum_{0\leq i\leq k} C_{i,k}\mu_i = \int
F(t)\sum_{0\leq i\leq k} C_{i,k}\phi_i(t) \dd t =$$ $$= \int
F(t)\psi_k(t)\dd t=\sum_{j=1}^N a_j \int f(t-x^j)\psi_k(t) \dd t =
\sum_{j=1}^N a_j\phi_k(x^j).$$ This completes the proof.
\end{proof}
\end{thm}
\medskip

According to Theorem \ref{thm:conv_dual}, in order to reduce the
reconstruction problem with the measurements (\ref{eq:measurments}) and for
signals of the a priori known form (\ref{eq:the_model}) to a solution of the
generalized Prony system we have to find an $f$ - convolution dual system
$\psi = \psi_k(t)$ to the measurements kernels $\phi$. In fact we need only
the coefficients $C_{i,k}$. Having these coefficients, we compute $M_k$ and
get system (\ref{eq:gnrlzd_prony}).

\smallskip

Solvability of system (\ref{eq:gnrlzd_prony}) and robustness of its solution
depend on the measurements kernels $\phi$. Specific examples are considered
below.

\subsection{Reconstruction from moments}\label{section_poly} We
are given a finite number of moments of a signal
$F(x)=\sum_{j=1}^Na_j f(x-x^j)$ as in
(\ref{eq:the_model}) in the form \be m_k=\int F(t)t^k \dd t,
\ k=(k_1,\dots,k_d)\geq (0,\dots,0).\ee So here
$\phi_k(x)=x^{k_1}_1 \cdots x^{k_d}_d$ for each multi-index $k$. We look
for the dual functions $\psi_k$ satisfying the convolution equation
\be\label{eq:conv_poly} \int f(t-x)\psi_k(t) \dd t=x^k \ee for each
multi-index $k$. To solve this equation we apply Fourier transform
to both sides of (\ref{eq:conv_poly}). Assuming that $\hat f(0)\neq 0$
and that $\hat f(\omega)$ has the derivatives of all the orders at $0$
we find (see \cite{Sar.Yom2})
that there is a unique solution to (\ref{eq:conv_poly}) provided by
\nc{\pa}{\partial}
\be\label{eq:dual_poly} \phi_k(x)=\sum_{l\leq k}C_{k,l}x^l, \ee
where
\[
C_{k,l}=\frac1{(\sqrt{2\pi})^d}{k\choose l} (-i)^{k+l} \left[
\left.\frac{\pa^{k-l}}{\pa\omega^{k-l}}
\right|_{\omega=0}\frac1{\hat f(\omega)}\right].
\]
The assumption $\hat f(0)\neq 0$ is essential in the construction of the
$f$-convolution dual system for the monomials as well as for other measurement
kernels (as well as in the study of the shifts of a given function in general).
The above calculation can be applied, with proper modifications, in more
general situations (see \cite{Sar.Yom2}).
On the other hand, the assumption of differentiability of $\hat f$ at zero
is not very restrictive, in particular, if we work with signals with finite
support.

\medskip

Returning to the moments reconstruction, we set the generalized
polynomial moments $M_k$ as
\be\label{eq:gnrlzd_moments_poly} M_k=\sum_{l\leq k}C_{k,l}m_l \ee
and obtain, as in Theorem \ref{thm:conv_dual}, the following
system of equations: \be\label{eaution_prony_poly}
\sum^N_{j=1}a_j(x^j)^k = \sum^N_{j=1}a_j (x^j_1)^{k_1}
\cdots (x^j_d)^{k_d}=M_k, \ k \geq (0,\dots,0).\ee This system is
called ``multidimensional Prony system". It appears in numerous
problems of theoretical and applied mathematics. In
Section \ref{sec:Prony.System} below
we recall one of the classical methods for its solution in
one-dimensional case, and describe, under certain assumptions, a
method for its solution in several dimensions. Local stability of the
solutions of the one-dimensional Prony system is discussed in
Section \ref{sec:prony}.

\subsection{Fourier case}\label{section_fourier}
The signal has the same form as in section \ref{section_poly}: 
$F(x)=\sum_{j=1}^Na_j f(x-x^j)$. The measurement kernels $\phi$ 
we now choose as the harmonics $\phi_k(x)=e^{ikx}$. So our measurements are
the Fourier coefficients $c_k(F)= \hat f(k)=\int F(x)e^{ikx} \dd x.$ 

\medskip

Let us assume now that $f$ satisfies the condition $\hat f(k)\ne 0$ for all
integer $k$. Then the Fourier harmonics turn out to form essentially 
$f$-self-dual system: indeed, taking $\psi_k(x)=\frac1{\hat f(k)}e^{ikx}$ we 
get immediately that $\int f(t-x)\psi_k(t) \dd t=\int f(t-x)\frac1{\hat f(k)}e^{ikt}\dd t=
\frac{\hat f(k)}{\hat f(k)}e^{ikx}=\phi_{k}(x).$ 

\medskip

According to our general scheme we put now $M_k=\frac1{\hat f(k)}c_k(F)$ and get
a system of the form 

\be\label{eq:prony_fourier} \frac1{\hat
f(k)}c_k(F)=M_k=\sum^N_{j=1}a_je^{ikx^j}=\sum^N_{j=1}a_j(\rho^j)^{k}, 
\ k\geq (0,\dots,0), \ee where $\rho^j=e^{ix^j}$. This is once more a 
multidimensional Prony system as in (\ref{eaution_prony_poly}), with the 
nodes on the complex unit circle.

\subsection{Further extensions}\label{subsection_extension}
The approach above can be extended in the following directions:
\begin{enumerate}
\item Reconstruction of signals built from several functions or with the 
addition of dilations also can be investigated (a perturbation approach 
where the dilations are approximately 1 is studied in \cite{Sar}). 
\item Further study of ``convolution duality" in \cite{Sar} provides a certain
extension of the class of signals and measurements allowing for a 
closed - form reconstruction.
\end{enumerate}

\section{Prony system}\label{sec:Prony.System}

\subsection{One-dimensional case}\label{sec:One.dim.Prony.System}

Let us start with the classical case of one-dimensional Prony system:

\be\label{prony.one.dim} \sum^N_{j=1}a_j(x^j)^k = m_k, \ x^j \in
{\mathbb R}, \ k=0,1,\dots. \ee This system appears in many branches of
mathematics (see \cite{lyubich2004sylvester}). There are several solution
methods available, for instance direct nonlinear minimization
(see e.g. \cite{Eck}), the polynomial realization (original Prony method,
\cite{prony1795essai}) or the state-space approach (\cite{rao1992mbp}).
Let us give a sketch of a method based on Pad\'e approximation techniques
which is rather close to the original Prony method. To simplify the presentation
we shall assume that we know a priori that all the nodes $x^j$ are pairwise
distinct. General case is treated similarly (\cite{Nik.Sor}, see also \cite{Yom1}).

Consider a ``moments
generating function" $I(z)=\sum_{k=0}^\infty m_k z^k, \ z\in {\mathbb C}.$
Assuming the equations (\ref{prony.one.dim}) are satisfied for all $k=0,1,\dots$
and summing up the geometric progressions we get

\be\label{eq:gener.fn.1} I(z)=\sum_{j=1}^N {a_j\over {1-x^jz}}.\ee So $I(z)$
is a rational function of degree $N$ tending to zero as $z\rightarrow \infty$.
The unknowns $a_j$ and ${1\over {x^j}}$ in (\ref{prony.one.dim}) are the poles
and the residues of $I(z)$, respectively.

Now in order to find $I(z)$ explicitly from the first $2N$ moments $m_0,m_1,\dots,m_{2N}$
we use the Pad\'e approximation approach (see \cite{Nik.Sor}): write
$I(z)$ as ${{P(z)}\over {Q(z)}}$ with polynomials
$P(z)=A_0+A_1z+\dots+A_{N-1}z^{N-1}$ and $Q(z)=B_0+B_1z+\dots+B_Nz^N$
of degrees $N-1$ and $N$, respectively.

Multiplying by $Q$ we have $I(z)Q(z)=P(z)$. Now equating the coefficients
on both sides we get the following system of linear equations:

\medskip

$m_0B_0=A_0$

\smallskip

$m_0B_1+m_1B_0=A_1$

\smallskip

.............................

\smallskip

$m_0B_{N-1}+m_1B_{N-2}+\dots+m_{N-1}B_0=A_{N-1}$

\smallskip

$m_0B_{N}+m_1B_{N-1}+\dots+m_{N-1}B_1+m_N B_0=0$

\smallskip

$m_1B_{N}+m_2B_{N-1}+\dots+m_{N}B_1+m_{N+1}B_0=0$

\smallskip

..............................

\medskip

The rest of the equations in this system are obtained by further shifts
of the indices of the moments, and so they form a Hankel-type matrix.

Now, being a rational function of degree $N$, $I(z)$ is uniquely defined
by its first $2N$ Taylor coefficients (the difference of two such functions
cannot vanish at zero with the order higher than $2N-1$). We conclude that
the linear system consisting of the first $2N$ equations as above is uniquely
solvable up to a common scaling of $P$ and $Q$ (of course, this fact follows
also form a general Pad\'e approximation theory - see \cite{Nik.Sor}).

\medskip

Now a solution procedure for the Prony system can be described as follows:

\medskip

1. Solve a linear system of the first $2N$ equations as above (with the
coefficients - the known moments $m_k$) to find the moments generating
function $I(t)$ in the form $I(z)={{P(z)}\over {Q(z)}}$.

\smallskip

2. Represent $I(z)$ in a standard way as the sum of elementary fractions
$I(z)=\sum_{j=1}^N {a_j\over {1-x^jz}}$. (Equivalently, find poles and residues
of $I(z)$). Besides algebraic operations, this requires just finding the roots
of the polynomial $Q(z)$. Then $(a_j, x^j), \ j=1,\dots,N$ form the unique
solution of the Prony system (\ref{prony.one.dim}).

\subsection{Multi-dimensional case}\label{sec:Mult.Prony.System}

Let us recall our multi-dimensional notations: 
$x=(x_1,\dots,x_d)\in {\mathbb R}^d$, \ $j$ \ is a scalar index, while 
$k=(k_1,\dots,k_d), \ i=(i_1,\dots,i_d)$ and $n=(n_1,\dots,n_d)$ are 
multi-indices. Partial ordering of multi-indices is given by 
$k \leq k' \Leftrightarrow k_p \leq k'_p, \ p=1,\dots,d.$

\smallskip

Let $x^j=(x^j_1,\dots,x^j_d), \ j=1,\dots,N.$ With the above notation, the
multidimensional Prony system has exactly the same form as in the one-dimensional
case:

\be\label{prony.mult.dim} \sum^N_{j=1}a_j(x^j)^k = m_k, \ x^j \in
{\mathbb R}^d, \ k\geq (0,\dots,0). \ee

Exactly as above we get that for $z=(z_1,\dots,z_d)\in {\mathbb C}^d$ 
the moments generating function $I(z)=\sum_{k\in \Nn^d}m_kz^k$ is a rational 
function of degree $Nd$ of the form

\be\label{eq:genrating_function} I(z)=
\sum_{j=1}^N a_j\prod_{l=1}^d\frac{1}{1-x^j_lz_l} \ . \ee Representing $I(z)$
as ${{P(z)}\over {Q(z)}}$ we get exactly in the same way as above an infinite 
system of linear equations for the coefficients of $P$ and $Q$, with a 
Hankel-type matrix formed by the moments $m_k$. By the same consideration as 
above, after we take enough equations in this system the solution is unique up
to a rescaling (see \cite{Bak,Nik.Sor} and references therein). 

\medskip

However, from this point the multi-dimensional situation becomes essentially
more complicated. While in dimension one $I(z)$ can be, essentially, any rational 
function of degree $N$ (naturally represented as the sum of elementary fractions), 
in several variables $I(z)$ has a very special form given by 
(\ref{eq:genrating_function}). This fact can be easily understood via counting the 
degrees of freedom: the number of unknowns in (\ref{prony.mult.dim}) is $N(d+1)$ 
while a rational function of degree $Nd$ in $d$ variables has much more coefficients than that, for $d>1$. It would be desirable to use as 
many equations from (\ref{prony.mult.dim}) as the number of unknowns, but the
method outlined above ignores a specific structure of $I(z)$ and requires as many
equations as in a Pad\'e reconstruction for a general rational function of degree
$Nd$.

\medskip

We treat this problem in \cite{Sar.Yom2} analyzing the structure of singularities
of $I(z)$ and on this base proposing a robust reconstruction algorithm. Let us
give here one simple special case of this algorithm, which separates the
variables in the problem.

\subsubsection{Separation of variables in the multi-dimensional Prony
system}\label{sec:Sep.Var.Prony.System}

Let us assume that we know {\it a priori} that the solution $(a_j,x^j), \ j=1,\dots,N$
of multi-dimensional Prony system (\ref{prony.mult.dim}) is such that all the
coordinates $x^j_l, \ l=1,\dots, d$ of the points
$x^j, \ j=1,\dots N$ are pairwise distinct. Moreover, we assume that
$a_{j_1}\ne a_{j_2}$ for $j_1\ne j_2$. Under these assumptions we proceed as follows:

\medskip

Consider ``partial moment generating functions"
$I_m(t), \ t\in {\mathbb C}, \ m=1,\dots,d$, defined by

\be I_m(t) = \sum_{r=1}^\infty m_{re_m}t^r,\ee where $e_m$ is a multi-index
with $(e_m)_j=0$ for $m\neq j$ and $1$ otherwise. We have the following simple
fact:

\begin{prop} $I_m(t)$ is a one-dimensional moments generating function of the form

\be\label{eq:part.gen.func} I_m(t)=\sum_{j=1}^Na_j\frac1{1-x^j_mt}.\ee It coincides
with the restriction of $I(z)$ to the $m$-th coordinate axis in ${\mathbb C}^d$.
\begin{proof} Let us evaluate
$I(z)$ along the $m$-th coordinate axis, that is on the line
$z=te_m$ with $e_m$ as above and $t\in {\mathbb C}$. We get
\[
I(te_m)=\sum_{j=1}^Na_j\prod_{l=1}^d\frac1{1-x^j_lt({e_m})_l}
=\sum_{j=1}^Na_j\frac1{1-x^j_mt}
\]
which is the moments generating function (\ref{eq:part.gen.func}). Now, to express
$I_m(t)$ through the multi-dimensional moments $m_k$ we notice that any monomial
$x^k$, with $k\ne pe_m$ vanishes identically on the $m$-th coordinate axis. Hence

\[
I(te_m)=\sum_{r=0}^\infty m_{r e_{m}}t^r.
\]
This shows that $I_m(t)\equiv I(te_m)$ and completes the proof of the proposition.
\end{proof}
\end{prop}
\medskip

Now applying the method described in Section (\ref{sec:One.dim.Prony.System})
above we find for each $m=1,\dots,d$ the coordinates $x^1_m,\dots,x^N_m$ and
(repeatedly) the coefficients $a_1,\dots,a_N$. It remains to arrange these
coordinates into the points $x^j=(x^j_1,\dots,x^j_d)$. This presents a certain
combinatorial problem, since Prony system (\ref{prony.mult.dim}) is invariant
under permutations of the index $j$. Under the assumptions above we proceed as
follows: for each $m=1,\dots,d$ we have obtained the (unordered) collection of
the pairs $(a_j, x^j_m), \ j=1,\dots,N$. By assumptions $a_{j_1}\ne a_{j_2}$
for $j_1\ne j_2$. Hence we can arrange in a unique way all the pairs
$(a_j, x^j_m), \ j=1,\dots,N, \ m=1,\dots,d$ into the string
$[(a_1,x^1_1),\dots,(a_1,x^1_d)],\dots,[(a_N,x^N_1),\dots,(a_N,x^N_d)].$ This
gives us the desired solution of multi-dimensional Prony system
(\ref{prony.mult.dim}).

\medskip

Notice that the assumption $a_{j_1}\ne a_{j_2}$ for $j_1\ne j_2$ is essential
here. Indeed, for $x_1\ne x_2$ and 
$x^1=(x_1,x_2), \ x^2=(x_2,x_1), \ \hat x^1=(x_1,x_1), \ \hat x^2=(x_2,x_2)$ 
we have $m_k =(x^1)^k+(x^2)^k \equiv (\hat x^1)^k+(\hat x^2)^k$ for $k$ on each
of the coordinate axes. So the (unique up to permutations of the index $j$) 
solution of the Prony system cannot be reconstructed from these moments only.

\medskip

Another remark is that the separation of variables as described above requires
knowledge of $2dN$ moments $m_k$ ($2N$ on each of the coordinate axes). This is
almost twice more than $N(d+1)$ unknowns. This number can be significantly reduced
in some cases. See \cite{Sar.Yom2} for further investigation in both of these 
directions.

\medskip

Stability estimates for the solution of one-dimensional Prony system (Section \ref{sec:prony}) can be
easily extended to the case considered in the present subsection. We do not give
here explicit statement of this result. 

\section{Reconstruction of piecewise $D$-finite functions from moments}\label{sec:dfinite}
In this section we present an overview of our previous findings on algebraic reconstruction of a certain general class of signals. See \cite{Bat} for the complete details.

Let \f{\fn:[a,b] \to \reals} consist of $\np+1$ ``pieces''
\f{\fn_0,\dotsc \fn_\np} with \f{\np \geq 0} jump points
\[
a=\jp_0 < \jp_1 \dotsc <\jp_{\np} < \jp_{\np+1}=b
\]
Furthermore, let $\fn$ satisfy on each continuity interval some
linear homogeneous differential equation with polynomial
coefficients: \f{\Op \fn_n \equiv 0, \; n=0,\dotsc,\np} where
\begin{equation}\label{eq:operator}
\mathfrak{D} = \sum_{j=0}^\theorder
\biggl(\sum_{i=0}^{\polymaxorder_j} \polycoeff_{i,j} x^i \biggr)
\frac{\dd^j}{\dd x^j} \quad (\polycoeff_{ij} \in \reals)
\end{equation}
Each \f{\fn_n} may be therefore written as a linear combination of
functions \f{\{u_i\}_{i=1}^\theorder} which are a basis for the
space \f{\nullsp{\Op} = \{f: \Op f \equiv 0\}}:
\begin{equation}\label{eq:linear-comb-nullsp}
\fn_n(x)=\sum_{i=1}^\theorder \alpha_{i,n} u_i(x), \quad n=0,1,\dotsc,\np
\end{equation}
We term such functions $\fn$ ``piecewise $D$-finite''. Many ``simple'' functions may be represented in this framework, such as polynomials, trigonometric polynomials and algebraic functions.

The sequence \f{\{\meas=\meas(\fn)\}} is given by the usual moments
\[
\meas(\fn) = \int_a^b x^k \fn(x) \dd x
\]

\begin{RP} Given \f{\theorder, \{\polymaxorder_i\}, \np,a,b} and
the moment sequence \f{\{m_k\}} of a piecewise $D$-finite function
$\fn$, reconstruct all the parameters
\f{\{\polycoeff_{i,j}\},\{\jp_i\}, \{\alpha_{i,n}\}}.
\end{RP}

The above problem can be solved as follows.
\begin{enumerate}
\item It turns out that the moment sequence of every piecewise $D$-finite function $\fn$ satisfies a linear recurrence relation, such that the coefficients of $\Op$ annihilating every piece of $\fn$ and the discontinuity locations $\{\jp_j\}$ can be recovered from the moments by solving a linear systems of equations plus a nonlinear step of polynomial roots finding. As we shall explain below, in many cases this step is equivalent to solving a certain generalized form of the previously mentioned Prony system \eqref{prony.one.dim}.
\item The function $\fn$ itself can be subsequently reconstructed by numerically calculating the basis for the space $\nullsp{\Op}$ and solving an additional linear system of equations to recover the coefficients $\alpha_{i,n}$.
\end{enumerate}

The above algorithm has been tested on reconstruction of
piecewise polynomials, piecewise sinusoids and rational functions - see \cite{Bat} for details. The results reported there were relatively accurate for low noise levels. In this paper we continue to explore the numerical stability of the method - see Sections \ref{sec:prony} and \ref{sec:numerical}.

Now let us show how the Prony system arises in the piecewise $D$-finite reconstruction method. Consider the distribution $\Op \fn$. It is of the form (see \cite[Theorem 2.12]{Bat})
\begin{equation}\label{eq:delta-der-model}
\Op \fn = \sum_{j=1}^{\np}\sum_{i=0}^{l_{j}-1}\jc_{i,j}\der{\delta}{i}\left(x-\jp_{j}\right)
\end{equation}
where $\jp_j \in \reals$ are the discontinuity locations, $\jc_{i,j} \in \reals$ are the associated ``jump magnitudes'' which depend on the values $\{\der{\fn}{i}(\jp_j^{\pm})\}$,  and $\delta(x)$ is the Dirac $\delta$-function. In particular, when $\fn$ is just a piecewise constant function, we have $\Op={\dd \over \dd x}$ and so
\begin{equation}\label{eq:delta-model}
\Op \fn = \sum_{j=1}^{\np}\jc_{j}\delta(x-\jp_{j})
\end{equation} where $\jc_j = \fn(\jp_j^{+})-\fn(\jp_j^{-})$. Let us now apply the moment transform to the equations \eqref{eq:delta-model} and \eqref{eq:delta-der-model}. We get, correspondingly,
\begin{eqnarray}
\meas\left(\Op \fn \right) &=&\sum_{j=1}^{\np}\jc_{j}\jp_{j}^{k}; \label{eq:prony}\\
\meas\left(\Op \fn \right) &=&\sum_{j=1}^{\np}\sum_{i=0}^{l_{j}-1}\jc_{i,j}k(k-1)\times\dots\times(k-i+1)\jp_{j}^{k-i}.\label{eq:gen-prony}
\end{eqnarray}
The system \eqref{eq:prony} is of course identical to the previously considered \eqref{prony.one.dim}. 
However, the following question arises: how are the numbers in the left-hand side of \eqref{eq:prony} and \eqref{eq:gen-prony} related to the known quantities $m_k\left(\fn\right)$? It turns out that the numbers $\meas\left(\Op \fn \right)$ are certain linear combinations of these known moments, with coefficients which are determined by $\Op$ in a well-defined way. The conclusion is as follows: \emph{if the operator $\Op$ which annihilates every piece of $\fn$ is known (but the other parameters are not\footnote{This assumption is realistic, for instance when reconstructing piecewise-polynomials or piecewise-sinusoids.}), then one can recover the discontinuity locations of $\fn$ by solving the Prony-like system \eqref{eq:gen-prony}.}
We term the system \eqref{eq:gen-prony} ``confluent Prony system''. It can be solved in a similar manner to the standard Prony system. A unique solution exists whenever all the $\jp_j$'s are distinct and the highest coefficients $\jc_{l_j-1,j}$ are nonzero. We provide the details in \cite{BatProny}. 

\section{Accuracy of solving Prony-like systems}\label{sec:prony}
In the preceding sections we have shown that several algebraic methods for nonlinear reconstruction can be reduced to solving certain types of systems of equations, the most basic one of which is the Prony system \eqref{prony.one.dim}, \eqref{eq:prony}. A crucial factor for the eventual applicability of the reconstruction methods is the sensitivity of solving these systems to measurement errors.

In this section we consider two such systems - \eqref{eq:prony} and \eqref{eq:gen-prony} and provide some theoretical results regarding the local sensitivity of their solution to noise. These results will hopefully enable further ``global" analysis.

We consider the following question: if the measurements in the left-hand side of \eqref{eq:prony} are known with error at most $\varepsilon$, how well can one recover the unknowns $\{\jc_j, \jp_j\}$? We regard this stability problem to be absolutely central in Algebraic Sampling. To our best knowledge, no general treatment of this problem exists, therefore we consider our results below to be a step in this direction.

For simplicity, let us assume that the number of equations in \eqref{eq:prony} equals the number of unknowns, which in this case is $P=2\np$. Let us further consider the ``measurement map'' \f{\fwm: \reals^P \to \reals^P} given by \eqref{eq:prony} (we call it the ``Prony map''):
\[
\fwm \bigl( \{ \jc_{j}\}, \{\jp_j\} \bigr) \isdef \{\meas\}_{k=0}^{P-1}
\]

Then, one possible answer (by no means a complete one) to the question posed above can be given in terms of the local Lipschitz constant of the ``solution map'' $\fwm^{-1}$, whenever this inverse is defined. We then obtain the following result.
\begin{thm}\label{thm:prony-accuracy}
Let $\{\meas\}_{k=0,\dots,P-1}$ be the exact unperturbed moments of the model \eqref{eq:delta-model}. Assume that all the $\jp_j$'s are distinct and also $\jc_j \neq 0$ for $j=1,\dots,\np$. Now let $\nn{\meas}$ be perturbations of the above moments such that $\max_k |\meas-\nn{\meas}| < \varepsilon$. Then, for sufficiently small $\varepsilon$, the perturbed Prony system has a unique solution which satisfies:
\begin{align*}
 |\nn{\jp}_j-\jp_j| &\leq C_1 \varepsilon |\jc_j|^{-1}\\
|\nn{\jc}_j - \jc_j| &\leq C_1 \varepsilon
\end{align*}
where $C_1$ is an explicit constant depending only on the geometry of \f{\jp_1,\dotsc,\jp_{\np}}.
\begin{proof}
Consider the Jacobian determinant of the map $\fwm$. It is easily  factorized as follows:
\[
\dd \fwm =
\underbrace{
\begin{bmatrix}
1 & 0 & \dots & 1 & 0\\
\jp_1 & 1 & \dots & \jp_{\np} & 1\\
\vdots\\
\jp_1^P & P \jp_1^{P-1} & \dots & \jp_{\np}^P & P \jp_{\np}^{P-1}
\end{bmatrix}}_{\isdef V \left( \jp_1,\dots,\jp_{\np} \right)}  \times \diag \{D_1,\dots,D_{\np} \}
\]
where $D_j$ is a $2\times 2$ block
\[
D_j \isdef
\begin{bmatrix}
1 & 0\\
0 & \jc_j
\end{bmatrix}
\]

The matrix $V=V \left( \jp_1,\dots,\jp_{\np} \right)$ is a special case of the so-called confluent Vandermonde matrix, well-known in numerical analysis (\cite{bjorck1973acv,gautschi1962iva, Gautschi1978445, kalman1984gvm}). In particular, the paper \cite[Theorem 3]{gautschi1962iva} contains the following estimate for the norm of $V^{-1}$:
\[
\|V^{-1}\|_{\infty} \leq \max_{1 \leq i \leq \np} b_{i} \prod_{j=1, j \neq i}^\np \biggl( \frac{1+|\jp_j|}{|\jp_i-\jp_j|} \biggr)^2
\]
where
\[
b_i\isdef \max \biggr(1+|\jp_i|,1+2(1+|\jp_i|) \sum_{j\neq i}\frac{1}{|\jp_j-\jp_i|}\biggl)
\]
Now since the $\jp_j$'s are distinct and $\jc_j\neq 0$, the Jacobian is non-singular and so in a sufficiently small neighbourhood of the exact solution, the map $\fwm$ is approximately linear. By the inverse function theorem
\[
\dd \fwm^{-1} = \diag \{D_1^{-1},\dots,D_{\np}^{-1} \} \times V^{-1}
\]
and so taking $C_1 \isdef \|V^{-1}\|_{\infty}$ completes the proof.
\end{proof}
\end{thm}

By a similar technique with slightly more involved calculations we obtain the following result for the system \eqref{eq:gen-prony}.
\begin{thm}
Let $\{\meas\}_{k=0,\dots,P-1}$ be the exact unperturbed moments of the model \eqref{eq:delta-der-model}, where $P=\sum_{j=1}^{\np} l_j+\np$. Assume that all the $\jp_j$'s are distinct and also $\jc_{l_j-1,j} \neq 0$ for $j=1,\dots,\np$. Now let $\nn{\meas}$ be perturbations of the above moments such that $\max_k |\meas-\nn{\meas}| < \varepsilon$. Then, for sufficiently small $\varepsilon$, the perturbed confluent Prony system has a unique solution which satisfies:
\begin{align*}
 |\nn{\jc}_{i,j}- \jc_{i,j}| &\leq \begin{cases}
                           C_2 \varepsilon & i=0\\ 
                           C_2 \varepsilon \biggl(1+\frac{|\jc_{i-1,j}|}{|\jc_{l_j-1,j}|}\biggr) & 1 \leq i \leq l_j-1
                          \end{cases}\\
|\nn{\jp}_j - \jp_j| &\leq C_2 \varepsilon \frac{1}{|\jc_{l_j-1,j}|}
\end{align*}
where $C_2$ is a constant depending only on the nodes \f{\jp_1,\dotsc,\jp_{\np}} and the multiplicities \f{l_1,\dotsc,l_{\np}}.
\begin{proof}[Proof outline]
As before, we obtain a factorization of the Jacobian determinant $\dd \fwm$ as a product of a confluent Vandermonde matrix $V(\jp_1,l_1,\dots,\jp_{\np},l_{\np})$ (defined in a similar manner but with each column having $l_j-2$ ``confluencies'') and the block diagonal matrix $D=\diag \{D_1,\dots, D_{\np} \}$ where
\[
D_j = \begin{bmatrix}
1 & 0 & \cdots & 0\\
0 & 1 & \cdots & \jc_{0,j}\\
\vdots & \vdots &\ \ddots & \vdots\\
0 & 0 & \cdots & \jc_{l_j-1,j}
\end{bmatrix}
\]
Inverting this expression and taking $C_2 = \|V^{-1}\|_{\infty}$ completes the proof.
\end{proof}
\end{thm}

An important generalization would be to consider the mappings
\[
\fwm_s : \bigl( \{ \jc_{i,j}\}, \{\jp_j\} \bigr) \isdef \{\meas\}_{k=s}^{P+s-1}
\]
and investigate the reconstruction error as $s\to\infty$. Such a formulation makes sense in the Fourier reconstruction problem (\cite{BatYomFourer, Eck, Ett.Sar.Yom}). This is a work in progress and we plan to present the results separately (\cite{BatProny}).

\section{Numerical experiments}\label{sec:numerical}
\begin{figure}[h]
\centering
\includegraphics[scale=0.6]{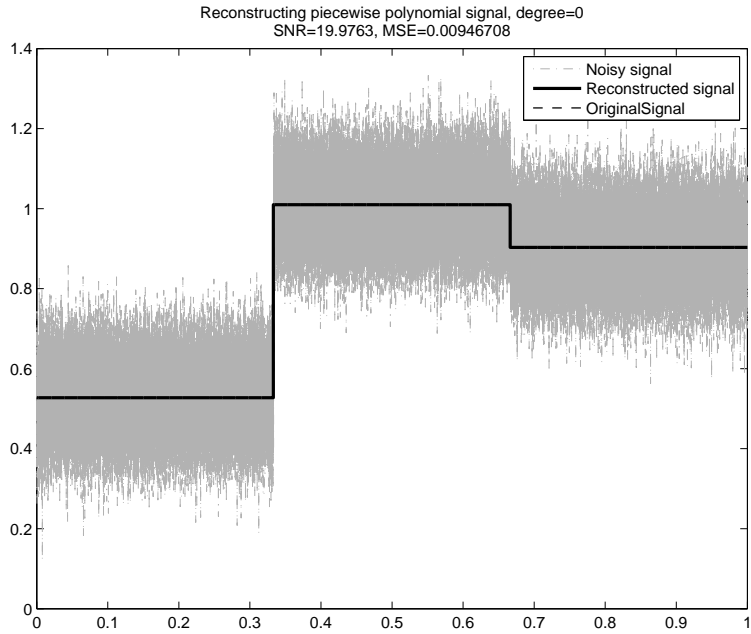}
\caption{Reconstructing a piecewise-constant signal with two jumps. The reconstructed and the original are very close.}
\label{fig:rectest}
\end{figure}

We have tested the piecewise $D$-finite reconstruction method on a simple case of a piecewise-constant signal. The implementation details are identical to those used in \cite[Appendix]{Bat}. As can be seen from Figure \ref{fig:rectest}, the reconstruction is accurate even in the presence of medium-level noise.

\begin{figure}[h]
\centering
\includegraphics[scale=0.6]{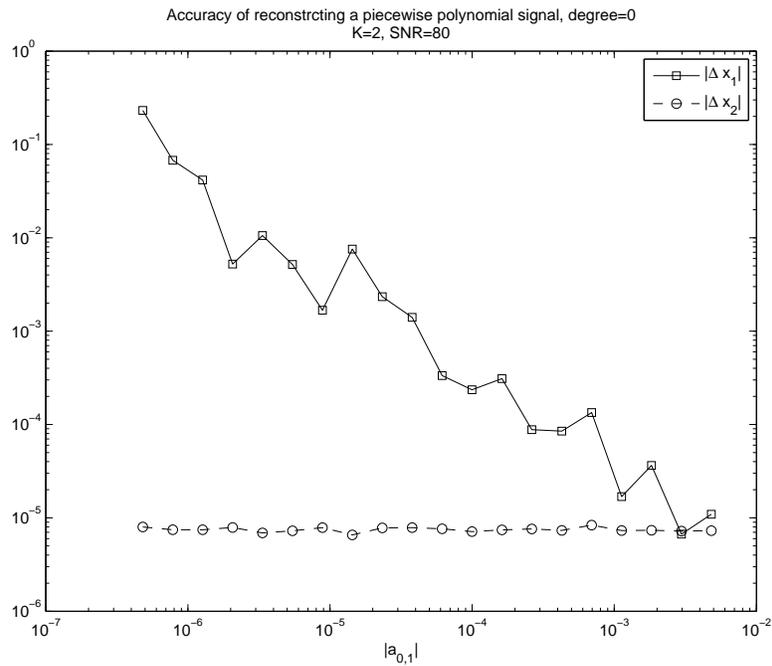}
\caption{Dependence of the reconstruction accuracy on the magnitude of the jump.}
\label{fig:highcoeff}

\end{figure}

We have already mentioned that for a piecewise constant function $\fun$, the distribution $\fun'$ is of the form \eqref{eq:delta-model}. Therefore, the piecewise $D$-finite  reconstruction problem for $\fun$ essentially reduces to solving the Prony system \eqref{eq:prony}, and so the local estimates of Theorem \ref{thm:prony-accuracy} should apply in the case of a very small noise level. This prediction is partially confirmed by the results of our second experiment, presented in Figure \ref{fig:highcoeff}. In particular, it can be seen that indeed $|\Delta \jp_j| \sim |\jc_j|^{-1}$, while the accuracy of other jump points $|\Delta \jp_i|$ does not depend on $|\jc_j|$ for $j\neq i$. 
While this is certainly an encouraging result, more investigation is clearly needed in order to fully understand the dependence of the error on all the parameters of the problem. Such an investigation should, in our opinion, concentrate on the global structure of the mapping $\fwm$.

%
%

\end{document}